\documentclass[11pt]{amsart}
\usepackage{latexsym,amsmath,amssymb}

\title[Lipschitz mappings and isometric embeddings]{{\protect{Sobolev mappings:
Lipschitz density is not an isometric invariant of the target}}}

\author{Piotr Haj\l{}asz}

\address{Department of Mathematics, University of Pittsburgh, 301
  Thackeray Hall, Pittsburgh, PA 15260, USA, {\tt hajlasz@pitt.edu}}

\thanks{The author was supported by NSF grant DMS-0900871.}

plus 6pt minus 12pt
\def\eps{\varepsilon}

\def\vi{\varphi}

\newtheorem{theorem}{Theorem}
\newtheorem{lemma}[theorem]{Lemma}
\newtheorem{corollary}[theorem]{Corollary}
\newtheorem{proposition}[theorem]{Proposition}

\def\lip{{\rm Lip\,}}

\theoremstyle{definition}

\newtheorem{definition}[theorem]{Definition}

\newcommand{\barint}{
\rule[.036in]{.12in}{.009in}\kern-.16in \displaystyle\int }

\newcommand{\barcal}{\mbox{$ \rule[.036in]{.11in}{.007in}\kern-.128in\int $}}

\def\eqn#1$$#2$${\begin{equation}\label#1#2\end{equation}}
\def\vint#1_#2{-\kern-#1pt\int_{#2}}

\newcommand{\bbbr}{\mathbb R}

\def\bbbc{{\mathchoice {\setbox0=\hbox{$\displaystyle\rm C$}\hbox{\hbox
to0pt{\kern0.4\wd0\vrule height0.9\ht0\hss}\box0}}
{\setbox0=\hbox{$\textstyle\rm C$}\hbox{\hbox
to0pt{\kern0.4\wd0\vrule height0.9\ht0\hss}\box0}}
{\setbox0=\hbox{$\scriptstyle\rm C$}\hbox{\hbox
to0pt{\kern0.4\wd0\vrule height0.9\ht0\hss}\box0}}
{\setbox0=\hbox{$\scriptscriptstyle\rm C$}\hbox{\hbox
to0pt{\kern0.4\wd0\vrule height0.9\ht0\hss}\box0}}}}

\def\bbbq{{\mathchoice {\setbox0=\hbox{$\displaystyle\rm Q$}\hbox{\raise
0.15\ht0\hbox to0pt{\kern0.4\wd0\vrule height0.8\ht0\hss}\box0}}
{\setbox0=\hbox{$\textstyle\rm Q$}\hbox{\raise 0.15\ht0\hbox
to0pt{\kern0.4\wd0\vrule height0.8\ht0\hss}\box0}}
{\setbox0=\hbox{$\scriptstyle\rm Q$}\hbox{\raise 0.15\ht0\hbox
to0pt{\kern0.4\wd0\vrule height0.7\ht0\hss}\box0}}
{\setbox0=\hbox{$\scriptscriptstyle\rm Q$}\hbox{\raise
0.15\ht0\hbox to0pt{\kern0.4\wd0\vrule height0.7\ht0\hss}\box0}}}}

\def\bbbz{{\mathchoice {\hbox{$\sf\textstyle Z\kern-0.4em Z$}}
{\hbox{$\sf\textstyle Z\kern-0.4em Z$}} {\hbox{$\sf\scriptstyle
Z\kern-0.3em Z$}} {\hbox{$\sf\scriptscriptstyle Z\kern-0.2em
Z$}}}}

\def\mvint_#1{\mathchoice
          {\mathop{\vrule width 6pt height 3 pt depth -2.5pt
                  \kern -8pt \intop}\nolimits_{\kern -3pt #1}}%
          {\mathop{\vrule width 5pt height 3 pt depth -2.6pt
                  \kern -6pt \intop}\nolimits_{#1}}%
          {\mathop{\vrule width 5pt height 3 pt depth -2.6pt
                  \kern -6pt \intop}\nolimits_{#1}}%
          {\mathop{\vrule width 5pt height 3 pt depth -2.6pt
                  \kern -6pt \intop}\nolimits_{#1}}}

\numberwithin{theorem}{section} \numberwithin{equation}{section}

\begin{document}

\subjclass[2000]{46E35}

\sloppy

\begin{abstract}
If $M$ is a compact smooth manifold and $X$ is a compact metric
space, the Sobolev space $W^{1,p}(M,X)$ is defined through an
isometric embedding of $X$ into a Banach space. We prove that the
answer to the question whether Lipschitz mappings $\lip(M,X)$ are
dense in $W^{1,p}(M,X)$ may depend on the isometric embedding of
the target.
\end{abstract}

\maketitle

\section{Introduction}

Sobolev mappings between manifolds $W^{1,p}(M,N)$ play a
fundamental role in geometric variational problems like the
theory of harmonic mappings. Eells and Lemaire \cite{eellsl},
asked whether smooth mappings are dense in $W^{1,p}(M,N)$ and it
turns out that the answer to that question depends on the topology
of the manifolds, see e.g. \cite{bethuel1}, \cite{bethuelz},
\cite{hajlasz2}, \cite{schoenu1}, \cite{schoenu2}, for early
results. Finally a necessary and sufficient condition was
discovered in Hang and Lin, \cite{hangl2}. The theory of Sobolev
mappings between manifolds has been extended to the case of
mappings into metric spaces. Research in that direction was
initiated in the work of Ambrosio, \cite{ambrosio}, Gromov and
Schoen, \cite{gromovs}, Korevaar and Schoen, \cite{korevaars},
Capogna and Lin \cite{capognal} and Reshetnyak \cite{reshetnyak}
just to
name a few. Finally the theory was even extended to the case of
Sobolev mappings between metric spaces, see Heinonen and Koskela
\cite{heinonenk} and Heinonen, Hoskela, Shanmugalingam and Tyson
\cite{HKST}. It is natural to inquire what would be suitable
generalizations of density results known in the case of mappings
between smooth manifolds to the case of a metric target. This
problem was explicitly formulated in the work of Heinonen,
Koskela, Shanumgalingam and Tyson \cite[Remark~6.9]{HKST} and some
partial results have been obtained in \cite{DHLT},
\cite{hajlaszGAFA} and \cite{hajlaszMathAnn}. For a more detailed
introduction to the subject, see the survey paper
\cite{hajlaszSobolev}.

In this paper we consider the Sobolev space $W^{1,p}(M,X)$
of mappings from a smooth compact Riemannian manifold (with or without boundary)
into a compact metric
space $X$. Every metric space admits an isometric embedding into a Banach
space; if $X$ is separable (in particular if $X$ is compact)
it can be isometrically embedded into $\ell^\infty$
(the Kuratowski embedding).
The space $\ell^\infty=(\ell^1)^*$ is dual to a separable Banach space.
Thus we may assume that
a compact space $X$ is isometrically
embedded into a Banach space $V$, $X\subset V$, where
$V=Y^*$ is dual to a separable Banach space $Y$.
In what follows we assume that every Banach space
into which we embed $X$ is a dual to separable Banach space.
The vector valued Sobolev
space $W^{1,p}(M,V)$ is a Banach space and we may define
$$
W^{1,p}(M,X)=
\left\{ f\in W^{1,p}(M,V):\, \mbox{$f(x)\in X$ a.e.}\right\}\, .
$$
If $X$ is compact, then all mappings into $X$ are bounded as mappings into $V$
and therefore integrable. The compactness assumption is to avoid problems
with integrability of the mapping.

The space $W^{1,p}(M,X)$ is equipped with a metric inherited from the norm of
$W^{1,p}(M,V)$ and we may inquire whether for a given metric space $X$,
Lipschitz mappings $\lip(M,X)$ are dense in $W^{1,p}(M,X)$.
This is exactly the problem that was formulated in \cite[Remark~6.9]{HKST}.
It turns out that Sobolev mappings $W^{1,p}(M,X)$ can be defined
in an intrinsic way independent of the isometric embedding, see Proposition~\ref{T5.5}.
However, the metric in $W^{1,p}(M,X)$ does depend on the embedding.
A simple example is provided in \cite[p. 438]{hajlaszGAFA}.
Thus if $\lambda:X\to V$ is an isometric embedding of $X$ into $V$,
used to define the metric in $W^{1,p}(M,X)$, we should rather denote the space
by $W^{1,p}_\lambda(M,X)$, but in practice, the subscript $\lambda$
is often omitted. Regarding the space of Lipschitz mappings $\lip(M,X)$
there is no need to use subscript $\lambda$.
The following question is very natural:

\noindent
{\sc Question.}
{\em Does the answer to the question about the density of Lipschitz mappings
in $W^{1,p}(M,X)$ depend on the isometric embedding of $X$ into a Banach space?}

This problem arose soon after the publication of \cite{HKST} and it was explicitly
formulated in \cite[Question~1]{hajlaszGAFA}.
The following result gives a partial answer to that problem.

\begin{theorem}
\label{T1}
Let $M$ be a compact Riemannian manifold, $X$ a compact metric space and $1\leq p<\infty$.
If for some isometric embedding $\lambda:X\to V$
of $X$ into a Banach space $V$ Lipschitz mappings $\lip(M,X)$
are dense in $W^{1,p}_\lambda(M,X)$ in the following strong sense:
for every $f\in W^{1,p}_\lambda(M,X)$ and every $\eps>0$ there is
$g\in\lip(M,X)$ such that
$$
|\{x:\, f(x)\neq g(x)\}|<\eps
\quad
\mbox{and}
\quad
\Vert f-g\Vert_{1,p}<\eps\, ,
$$
then Lipschitz mappings are dense in $W^{1,p}_\nu(M,X)$ for any other
isometric embedding $\nu:X\to W$ of $X$ into a Banach space $W$.
\end{theorem}
As explained above, we assume here that each of the Banach spaces
$V$ and $W$ is  dual to a separable Banach space.
This result
is a version of \cite[Theorem~4]{hajlaszMathAnn}.
The density in the strong sense seems a typical property.
The following result was proved in \cite[Lemma~13]{hajlaszMathAnn}.
\begin{proposition}
\label{T1.5}
If a Banach space $V$ is dual to a separable Banach space and
$f\in W^{1,p}(M,V)$, then for every $\eps>0$ there is $g\in\lip(M,V)$
such that $|\{x:\, f(x)\neq g(x)\}|<\eps$ and $\Vert f-g\Vert_{1,p}<\eps$.
\end{proposition}
In view of the two results
one could naturally expect that the answer to the above question should be that
the density is not depend on the choice of the isometric embedding and that is
a completely incorrect intuition. Indeed, the following theorem,
the main result of the paper,
provides an example that gives a different answer.
\begin{theorem}
\label{T2}
There is a compact and connected set
$X\subset \bbbr^{n+2}$ such that
Lipschitz mappings $\lip(S^n,X)$ are dense in $W^{1,n}(S^n,X)$, while if
$\kappa:X\to\ell^\infty$ is the Kuratowski embedding, then
Lipschitz mappings $\lip(S^n,X)$ are not dense in $W^{1,n}_\kappa(S^n,X)$.
\end{theorem}
Since $X$ is a subset of $\bbbr^{n+2}$ in the first definition of
$W^{1,n}(S^n,X)$ we just consider the identity embedding of $X$ into
$\bbbr^{n+2}$ (which is a Banach space) and we avoid the subscript ${\rm id}$.

The paper is organized as follows.
In Section~\ref{SM} we review the theory of Sobolev mappings into metric spaces.
The main result of the section is Corollary~\ref{T7}
which is one of the main tools in the proof of Theorem~\ref{T2}.
In Section~\ref{Tejeden} we prove Theorem~\ref{T1}
and Section~\ref{Proof} is devoted to the proof of Theorem~\ref{T2}.

{\bf Acknowledgements.} The author acknowledges the kind hospitality
of CRM at the Universitat Aut\`onoma de Barcelona, where the research
was partially carried out.

\section{Sobolev mappings into metric spaces}
\label{SM}

In this section we briefly discuss the construction of the Sobolev
space of mappings from a manifold into a metric space.
We follow the presentation given in \cite{hajlaszt}.
We refer the reader to that paper for detailed proofs and references.
For simplicity  we formulate most of the definitions and results for
open subsets $\Omega\subset\bbbr^n$ instead of compact manifolds, but
the generalization to the manifold case is straightforward via the use
of local coordinate systems.
Most of the results in this section are known. The new results are
Lemma~\ref{T6} and Corollary~\ref{T7}.

If $V$ is any Banach space
(not necessarily dual to a separable Banach space)
and $A\subset\bbbr^n$ is (Lebesgue)
measurable, we say that $f\in L^p(A,V)$ if
\begin{enumerate}
\item[(1)] $f$ is {\em essentially separably valued}: $f(A\setminus
Z)$ is a separable subset of $V$ for some set $Z$ of Lebesgue measure zero,
\item[(2)] $f$ is {\em weakly measurable}: for every $v^*\in V^*$ with
  $\Vert v^*\Vert\leq 1$, $\langle v^*,f\rangle$ is measurable,
\item[(3)] $\Vert f\Vert\in L^p(A)$.
\end{enumerate}
If $f\in L^1(A,V)$ we define the integral
$$
\int_A f(x)\, dx\in V
$$
in the Bochner sense, see \cite[Chapter~5, Sections~4-5]{yosida},
\cite{diestelu}.
The Bochner integral has two important properties:
For every $v\in V^*$
$$
\left \langle v^*, \int_A f(x) \, dx \right \rangle = \int_A \langle
v^*, f(x) \rangle \, dx
$$
and
$$
\left\Vert \int_A f(x) \, dx \right\Vert \le \int_A ||f(x)|| \, dx.
$$
\begin{definition}
Let $\Omega\subset\bbbr^n$ be an open set and $V$ any
Banach space (not necessarily dual). The Sobolev space
$W^{1,p}(\Omega,V)$, $1\leq p<\infty$, is defined as the class of
all functions  $f\in L^p(\Omega,V)$ such that for $i=1,2,\ldots,n$
there is $f_i\in L^p(\Omega,V)$ such that for every $\vi\in
C_0^\infty(\Omega)$
$$
\int_\Omega\frac{\partial\vi}{\partial x_i}\, f = - \int_\Omega \vi f_i\, ,
$$
where the integrals are taken in the sense of Bochner (note that the
integrands are supported on compact subsets of $\Omega$). We denote
$f_i=\partial f/\partial x_i$ and call these functions {\em weak
  partial derivatives} of $f$. We also write $\nabla f =(\partial
f/\partial x_1,\ldots, \partial f/\partial x_n)$ and
\begin{equation}\label{gradient}
|\nabla f|=
\left( \sum_{i=1}^n \left\Vert\frac{\partial f}{\partial x_i}\right\Vert^2\right)^{1/2}\, .
\end{equation}
Sometimes we will write $|\nabla f|_V$ to emphasize the Banach space
with respect to which we compute the length of the gradient.
The space $W^{1,p}(\Omega,V)$ is equipped with the norm
$$
\Vert f\Vert_{1,p} = \left(\int_\Omega \Vert f\Vert^p\right)^{1/p} +
\left(\int_\Omega |\nabla f|^p\right)^{1/p}\, .
$$
It is easy to prove that $W^{1,p}(\Omega,V)$ is a Banach space.
\end{definition}

It easily follows from the definition (see \cite[Proposition~2.3]{hajlaszt})
that for every $v^*\in V^*$ with
$\Vert v^*\Vert\leq 1$, we have
$\langle v^*,f\rangle\in W^{1,p}(\Omega)$ and
\begin{equation}
\label{res}
|\nabla \langle v^*,f\rangle|\leq |\nabla f|
\quad
\mbox{a.e.}
\end{equation}
Observe that $v^*:V\to\bbbr$ is a $1$-Lipschitz function
and it turns out that under the additional assumption that $V$ is dual to a separable
Banach space (\ref{res}) holds with $v^*$ replaced by any
$1$-Lipschitz function. Moreover $|\nabla f|$ is, in a certain sense,
the best lower bound for $|\nabla \langle v^*,f\rangle|$. Namely we have.
\begin{proposition}
\label{T2.5}
Let $\Omega\subset\bbbr^n$ be open, let $V$ be dual to a separable
Banach space and let $1\leq p<\infty$. Then
for $f\in W^{1,p}(\Omega,V)$ we have
\begin{enumerate}
\item If $0\leq g\in L^p(\Omega)$ is such that for every $v^*\in V^*$
with $\Vert v^*\Vert\leq 1$ we have
$$
|\nabla \langle v^*,f\rangle|\leq g
\quad
\mbox{a.e.}
$$
then
$$
|\nabla f|\leq C g
\quad
\mbox{a.e.}
$$
\item For every $1$-Lipschitz function $\vi:V\to\bbbr$, we have
$\vi\circ f\in W^{1,p}(\Omega)$ and
$$
|\nabla (\vi\circ f)|\leq |\nabla f|.
$$
\end{enumerate}
\end{proposition}
The first part of this result follows from Theorem~2.14 and Lemma~2.13 in
\cite{hajlaszt}. The second part is a consequence of the estimate (\ref{res})
and the proof of Proposition~2.16 in \cite{hajlaszt}.

The proof of the above result is based on the characterization of Sobolev mappings
by absolute continuity on lines.
For the following result, see \cite[Lemma~2.8 and Lemma~2.13]{hajlaszt}.
\begin{lemma}
\label{T4}
Let $V=Y^*$ be dual to a separable Banach space $Y$. If $f:[a,b]\to V$
is absolutely continuous, then
the limit
$$
g(x):=\lim_{h\to 0}\left\Vert\frac{f(x+h)-f(x)}{h}\right\Vert
$$
exists a.e. and $g\in L^1([a,b])$. Moreover
for a.e. $x\in (a,b)$ there is a vector
$f'(x)\in V$ such that $\Vert f'(x)\Vert\leq g(x)$
and
$$
\left\langle v^*, \frac{f(x+h)-f(x)}{h}\right\rangle \to
\langle v^*,f'(x)\rangle
\qquad \mbox{as $h\to 0$}
$$
for all $v^*\in Y$.
We call $f'(x)$ the $w^*$-derivative of $f$ at $x$.
\end{lemma}
The lemma leads to the following characterization of the Sobolev space,
see \cite[Lemma~2.12, Lemma~2.13 and Theorem~2.14]{hajlaszt}.
\begin{proposition}
\label{T5}
Let $\Omega\subset\bbbr^n$ be an open set,
let $V=Y^*$ be dual to a separable Banach space $Y$ and let $1\leq p<\infty$.
Let $f\in L^p(\Omega,V)$.
Then $f\in W^{1,p}(\Omega,V)$ if and only if $f$ is absolutely continuous on
compact intervals in $\ell\cap\Omega$ for almost all lines $\ell$ parallel to
coordinate axes (possibly after being redefined on a set of measure zero)
and if $w^*$-partial derivatives of $f$ belong to  $L^p(\Omega,V)$.
Moreover the $w^*$-partial derivatives are equal to weak derivatives of $f$.
\end{proposition}
Now we are ready to define the Sobolev space of mappings with values into a
compact metric space $X$. If $\lambda:X\to V$ is an isometric embedding
of $X$ into a Banach space $V$ which is dual to a separable
Banach space, then we define
$$
W^{1,p}(\Omega,X)=W^{1,p}_\lambda(\Omega,X)=
\{f\in W^{1,p}(\Omega,V):\, f(x)\in\lambda(X)\ \mbox{a.e.}\}.
$$
Every compact (or even separable) metric space $(X,d)$ admits an isometric
embedding into $\ell^\infty$. Indeed, if
$\{x_i\}_{i=1}^\infty$ is a dense subset of $X$ and $x_0\in X$ any point, then
one can easily show that the mapping
$$
\kappa:X\to\ell^\infty,
\qquad
\kappa(x)=\big(d(x,x_i)-d(x_0,x_i)\big)_{i=1}^\infty
$$
is the isometric embedding. It is the well known {\em Kuratowski embedding}.
Therefore the Sobolev space $W^{1,p}(\Omega,X)$, can be defined for any compact
metric space $X$, because we can always use the Kuratowski embedding. Moreover observe that
$\ell^\infty = (\ell^1)^*$, so the space $\ell^\infty$ is dual to a separable
Banach space.

It turns out that the Sobolev space of mappings into $X$ can be defined in an
intrinsic way independent of the choice of the embedding $\lambda$.
For the following result see Definition~1.2 and argument on p. 698 in \cite{hajlaszt}.
\begin{proposition}
\label{T5.5}
Let $\Omega\subset\bbbr^n$ be a bounded open set,
$X$ a compact metric space and $1\leq p<\infty$.
Then $f\in W^{1,p}(\Omega,X)$ if and only if there is a nonnegative function
$g\in L^p(\Omega)$ such that for every Lipschitz function
$\vi:X\to\bbbr$, $\vi\circ f\in W^{1,p}(\Omega)$ and
$|\nabla(\vi\circ f)|\leq\lip(\vi) g$ a.e.
Here $\lip(\vi)$ stands for the Lipschitz constant of $\vi$.
\end{proposition}

Let $f=(f^i)_{i=1}^\infty\in W^{1,p}(\Omega,\ell^\infty)$.
Then $f$ is absolutely continuous on almost all lines parallel to coordinate axes.
Since the weak partial derivatives are equal to $w^*$-partial derivatives
we can compute them with help of Lemma~\ref{T4}.

Let $\delta_i\in (\ell^\infty)^*$ be defined as $i$th coordinate, i.e.
$$
\delta_i(z_1,z_2,\ldots)=z_i,
$$
Hence $f^i=\delta_i\circ f$ is absolutely continuous on almost all lines
as a composition of $f$ with the Lipschitz function $\delta_i$.
Thus the coordinate functions $f^i$ belong to $W^{1,p}(\Omega)$.

For $k=1,2,\ldots,n$, Lemma~\ref{T4} gives a formula for weak partial derivatives
$$
\frac{\partial f}{\partial x_k}=
\left(\left(\frac{\partial f}{\partial x_k}\right)^i\right)_{i=1}^\infty
\in L^p(\Omega,\ell^\infty)\, .
$$
Indeed,
$$
\frac{f^i(x+he_k)-f^i(x)}{h}=
\left\langle \delta_i, \frac{f(x+he_k)-f(x)}{h}\right\rangle \to
\left\langle \delta_i, \frac{\partial f}{\partial x_k}(x)\right\rangle=
\left(\frac{\partial f}{\partial x_k}\right)^i
$$
a.e., that is
$$
\left(\frac{\partial f}{\partial x_k}(x)\right)^i =
\frac{\partial f^i}{\partial x_k}(x)
\quad
\mbox{a.e.}
$$
The next two results are new.
\begin{lemma}
\label{T6}
Let $f,g:[a,b]\to\bbbr^N$ be absolutely continuous and let
$\kappa:\bbbr^N\to\ell^\infty$ be the Kuratowski embedding.
Then $\bar{f}=\kappa\circ f$ and $\bar{g}=\kappa\circ g$
are absolutely continuous functions with values into
$\ell^\infty$ and the $w^*$-derivative $(\bar{f}-\bar{g})':[a,b]\to\ell^\infty$
satisfies
$$
\Vert (\bar{f}-\bar{g})'(t)\Vert_\infty \geq
\max\{ |f'(t)|, |g'(t)|\} \geq
\frac{1}{2}\left( |f'(t)|+|g'(t)|\right)
$$
for almost every $t\in [a,b]$ such that $f(t)\neq g(t)$.
\end{lemma}
{\em Proof.}
Let $\{ x_i\}_{i=1}^\infty\subset\bbbr^N$
be a dense subset, $x_0\in\bbbr^N$ and let
$\kappa:\bbbr^N\to\ell^\infty$ be the Kuratowski embedding.
It is easy to see that
$$
(\bar{f}-\bar{g})(s)=
\left(|f(s)-x_i|-|g(s)-x_i|\right)_{i=1}^\infty
$$
and
$$
(\bar{f}-\bar{g})'(s)=
\left(|f(s)-x_i|'-|g(s)-x_i|'\right)_{i=1}^\infty
$$
for a.e. $s\in [a,b]$.
Fix $t\in [a,b]$ such that $f(t)\neq g(t)$ and both $f$ and $g$ are
differentiable at $t$.
If $f'(t)=g'(t)=0$, the inequality is obvious. Assume, then that one of the
derivatives is non zero, say $f'(t)\neq 0$.
Let $\delta>0$.
If we choose $x_i$ to be very close to $f(t)+f'(t)\delta$, then the function
$s\mapsto |f(s)-x_i|$ is decreasing near $s=t$ at the rate very close to
$|f'(t)|$, because the point $x_i$ is nearly exactly in the direction in which
$f(s)$ is going.
More precisely
$$
|f(t)-x_i|'=\frac{f'(t)\cdot (f(t)-x_i)}{|f(t)-x_i|}=|f'(t)|\cos\theta,
$$
where $\theta$ is the angle between the vectors $f'(t)$ and
$f(t)-x_i$. If $x_i$ is very close to $f(t)+f'(t)\delta$,
$\theta$ is very close to $\pi$.
Hence given $\eps>0$ and $\delta>0$ we can find $x_i$
so close to $f(t)+f'(t)\delta$ that
$$
|f(t)-x_i|'\leq -|f'(t)|+\eps.
$$
Choosing $x_{i'}$ close to $f(t)-f'(t)\delta$ we make the function
$s\mapsto |f(s)-x_{i'}|$ increasing at the rate very close to $|f'(t)|$,
so given $\eps>0$ and  $\delta>0$ we can find $x_{i'}$ so close to
$f(t)-f'(t)\delta$ that
$$
|f(t)-x_{i'}|'\geq |f'(t)|-\eps.
$$
Since $|f(t)-g(t)|>0$ (remember that we assume that $f(t)\neq g(t)$),
taking $\delta>0$ sufficiently small we can make the points $x_i$ and
$x_{i'}$ so close to $f(t)$ that
$$
\left| |g(t)-x_i|'-|g(t)-x_{i'}|'\right|<\eps.
$$
Hence either
$$
\left| |f(t)-x_i|'-|g(t)-x_i|'\right|\geq |f'(t)|-2\eps
$$
or
$$
\left| |f(t)-x_{i'}|'-|g(t)-x_{i'}|'\right|\geq |f'(t)|-2\eps
$$
by the triangle inequality. Thus
\begin{equation}
\label{wewe}
\Vert (\bar{f}-\bar{g})'(t)\Vert_\infty =
\sup_j \left| |f(t)-x_j|'-|g(t)-x_j|'\right| \geq
|f'(t)|-2\eps.
\end{equation}
Since the inequality is true for any $\eps>0$ we have
$$
\Vert (\bar{f}-\bar{g})'(t)\Vert_\infty\geq |f'(t)|.
$$
If $g'(t)=0$ the lemma follows. If $g'(t)\neq 0$ we can repeat the above argument
with $f$ replaced by $g$ and obtain the estimate
$$
\Vert (\bar{f}-\bar{g})'(t)\Vert_\infty\geq |g'(t)|.
$$
The two estimates combined together prove the lemma.
\hfill $\Box$

\begin{corollary}
\label{T7}
Let $\Omega\subset\bbbr^n$ be bounded.
Let $f,g\in W^{1,p}(\Omega,\bbbr^N)$ and let
$\kappa:\bbbr^N\to\ell^\infty$ be the Kuratowski embedding.
Then $\bar{f}=\kappa\circ f, \bar{g}=\kappa\circ g\in W^{1,p}(\Omega,\ell^\infty)$ and
$$
|\nabla (\bar{f}-\bar{g})|\geq
C(n)\left( |\nabla f|+|\nabla g|\right)\chi_{\{ f\neq g\}}
\quad
\mbox{a.e.}
$$
\end{corollary}
{\em Proof.}
The result follows immediately from the fact that $f$ and $g$ are absolutely continuous
almost all lines parallel to the coordinate directions, from Lemma~\ref{T4}
and from the definition
$$
|\nabla(\bar{f}-\bar{g})|=
\left(\sum_{k=1}^n
\left\Vert\frac{\partial (\bar{f}-\bar{g})}{\partial x_k}\right\Vert_\infty^2\right)^{1/2}\, .
$$
The proof is complete.
\hfill $\Box$

\section{Proof of Theorem~\ref{T1}}
\label{Tejeden}

The following result is well known in the case of real valued Sobolev functions,
but since the vector valued case is more delicate, we provide a short proof.
\begin{lemma}
\label{T7.5}
Let $V=Y^*$ be dual to a separable Banach space, let
$\Omega\subset\bbbr^n$ be an open set and let $1\leq p<\infty$. If
the functions $f_1,f_2\in W^{1,p}(\Omega,V)$ are equal on a measurable
set $E$, then $\nabla f_1=\nabla f_2$ a.e. on $E$.
\end{lemma}
{\em Proof.} Let $f=f_1-f_2$. Then $f=0$ on $E$ and we need to prove that
$\nabla f=0$ a.e. on $E$. Let $v^*\in Y$. Then
$\langle v^*,f\rangle\in W^{1,p}(\Omega)$. Since
$\langle v^*,f\rangle=0$ on $E$ it is well known that
$\nabla\langle v^*,f\rangle(x)=0$ for all $x\in E\setminus Z_{v^*}$ for some set
$Z_{v^*}$ of Lebesgue measure zero.
Let $D\subset Y$ be a countable dense set.
The set $Z=\bigcup_{v^*\in D} Z_{v^*}$ has measure zero and
$\nabla \langle v^*,f\rangle(x)=0$ for all $x\in E\setminus Z$
and all $v^*\in D$.
The weak partial derivatives of $f$ are equal to $w^*$-partial derivatives, see
Proposition~\ref{T5}.
Let $f$ be absolutely continuous on $\ell\cap\Omega$, where $\ell$ is parallel
to the $k$th coordinate axis. Then by Lemma~\ref{T4} for $v^*\in Y$
$$
\left\langle v^*,\frac{f(x+he_k)-f(x)}{h}\right\rangle
\to
\left\langle v^*,\frac{\partial f}{\partial x_k}(x)\right\rangle
\quad
\mbox{for a.e. $x\in\ell\cap\Omega$}.
$$
On the other hand the limit equals zero if $x\in\ell\cap (E\setminus Z)$ and $v^*\in D$,
so
$$
\left\langle v^*,\frac{\partial f}{\partial x_k}(x)\right\rangle = 0
\quad
\mbox{for $v^*\in D$ and $x\in\ell\cap(E\setminus Z)$}.
$$
Hence $\partial f/\partial x_k(x)=0$ for a.e. $x\in\ell\cap E$
by density of $D$ in $Y$.
\hfill $\Box$

Let $\nu:X\to W$ be an isometric embedding and
$f\in W^{1,p}_\nu(M,X)$. Then
$\bar{f}=\lambda\circ\nu^{-1}\circ f\in W^{1,p}_\lambda(M,X)$, because by
Proposition~\ref{T5.5}, the Sobolev space of mappings into $X$ can be defined
independently of the isometric embedding.
It follows from the assumptions of the theorem that
there is a sequence of Lipschitz mappings
$\bar{g}_k\in \lip(M,\lambda(X))$ such that
$$
|\{ x:\, \bar{f}(x)\neq \bar{g}_k(x)\}|\to 0
\quad
\mbox{and}
\quad
\Vert \bar{f}-\bar{g}_k\Vert_{1,p}\to 0
\quad
\mbox{as $k\to\infty$.}
$$
Then $g_k=\nu\circ\lambda^{-1}\circ \bar{g}_k\in \lip(M,\nu(X)))$ and
$$
|\{ x:\, f(x)\neq g_k(x)\}|=
|\{x:\, \bar{f}(x)\neq \bar{g}_k(x)\}|\to 0
\quad
\mbox{as $k\to\infty$.}
$$
In particular $g_k\to f$ in $L^p(M,W)$, because $X$ is a bounded subset of $W$.
It remains to estimate the gradients.

It follows from Lemma~\ref{T7.5} that $\nabla f=\nabla g_k$ a.e. on the set
where $f=g_k$, so
\begin{eqnarray*}
\lefteqn{\left(\int_M |\nabla f-\nabla g_k|^p_W\, dx\right)^{1/p}
 =
\left(\int_{\{f\neq g_k\}} |\nabla f-\nabla g_k|^p_W\, dx\right)^{1/p}}\\
& \leq &
\left(\int_{\{f\neq g_k\}} |\nabla f|^p_W\, dx\right)^{1/p} +
\left(\int_{\{f\neq g_k\}} |\nabla g_k|^p_W\, dx\right)^{1/p}\, .
\end{eqnarray*}
The first integral on the right hand side converges to zero, because
$|\{ f\neq g_k\}|\to 0$ and we need to show that the second integral
converges to zero as well.
For $w^*\in W^*$ with $\Vert w^*\Vert\leq 1$ the function
$$
v\mapsto \langle w^*,\nu(\lambda^{-1}(v))\rangle
$$
is $1$-Lipschitz continuous on $\lambda(X)\subset V$
and hence it extends to a $1$-Lipschitz continuous function
$\vi:V\to\bbbr$ (McShane extension). Since
$$
\langle w^*,g_k(x)\rangle =
(\vi\circ\bar{g}_k)(x),
\quad
x\in M
$$
Proposition~\ref{T2.5} gives
$$
|\nabla\langle w^*,g_k\rangle| =
|\nabla(\vi\circ \bar{g}_k)|\leq
|\nabla \bar{g}_k|_V
\quad
\mbox{a.e.}
$$
Then another application of Proposition~\ref{T2.5} yields
$$
|\nabla g_k|_W\leq C |\nabla \bar{g}_k|_V
\quad
\mbox{a.e.}
$$
Hence
\begin{eqnarray*}
\lefteqn{\left(\int_{\{f\neq g_k\}} |\nabla g_k|^p_W\, dx\right)^{1/p}
 \leq
C \left(\int_{\{f\neq g_k\}} |\nabla \bar{g}_k|^p_V\, dx\right)^{1/p}}\\
& \leq &
C\left(\left(\int_{\{f\neq g_k\}} |\nabla \bar{f}-\nabla \bar{g}_k|^p_V\, dx\right)^{1/p} +
\left(\int_{\{f\neq g_k\}} |\nabla \bar{f}|^p_V\, dx\right)^{1/p}\right)
\to 0
\end{eqnarray*}
as $k\to\infty$. The proof is complete.
\hfill $\Box$

\section{Proof of the main Theorem~\ref{T2}}
\label{Proof}

We begin with the construction of the set $X$. Actually this is exactly
the same set as in \cite{hajlaszGAFA} where it was used to provide
a counterexample to a different question.
In the description below we will try to emphasize the geometric nature of the
construction and we will avoid all technical details. Actually the details
of the construction and proofs are quite involved and we refer
the reader to \cite{hajlaszGAFA} for a detailed exposition.

In the first step we construct a continuous function $\gamma\in W^{1,n}$ on $\bbbr^n$
with compact support contained in $B(0,1)$, $\gamma(0)=0$.
The function is actually $C^\infty$ smooth except at the origin, where it has
accumulating oscillations with height gradually vanishing at $0$.

Next we replace a subset of $S^n$ diffeomorphic to $B(0,1)$ with the graph of
$\gamma$. The resulting space
denoted by $S_\infty$
is homeomorphic to $S^n$, and actually diffeomorphic
everywhere but at one point.
$S_\infty$ is constructed as a subset of $\bbbr^{n+1}$.
Since $\gamma$ belongs to $W^{1,n}$ there is a $W^{1,n}$ homeomorphism $f$ of
$S^n$ onto $S_\infty$.

It turns out that Lipschitz mappings
$\lip(S^n,S_\infty)$ are not dense in $W^{1,n}(S^n,S_\infty)$.
Indeed, homotopy properties of $W^{1,n}$ mappings  imply that if
a Lipschitz mapping $g\in \lip(S^n,S_\infty)$ is sufficiently close to
$f$, in the Sobolev norm, then $g$ must be a surjective mapping.
On the other hand the
oscillations of $\gamma$ are so frequent that there is no Lipschitz surjection
$g:S^n\to S_\infty$.

Thus there is $\eps>0$ such that
\begin{equation}
\label{far}
\mbox{$\Vert f-g\Vert_{1,n}>\eps$ for all $g\in\lip(S^n,S_\infty)$.}
\end{equation}
The function $\gamma$ is defined as a series
$$
\gamma=\sum_{i=1}^\infty \eta_i,
$$
where $\eta_k$ are a smooth, compactly supported bump functions.
Let $\widetilde{S}_k$ be the manifold obtained from $S^n$ be replacing
a subset diffeomorhphic to $B(0,1)$ with the graph of
$$
\gamma_k=\sum_{i=1}^k \eta_i.
$$
Clearly $\widetilde{S}_k\subset\bbbr^{n+1}$ is a smooth manifold
diffeomorphic to $S^n$ and the sequence
$\widetilde{S}_k$ converges in some sense to $S_\infty$.
Each set $\widetilde{S}_k$ and $S_\infty$ is a subset of
$\bbbr^{n+1}\subset\bbbr^{n+2}=\bbbr^{n+1}\times\bbbr$.
Let $S_k$ be the translation of $\widetilde{S}_k$ by the vector
$\langle 0,\ldots,0,2^{-k}\rangle$ in $\bbbr^{n+2}$
and now we define
$$
\widetilde{X}=S_\infty\cup\bigcup_{k=1}^\infty S_k.
$$
Thus the set $\widetilde{X}$ consists of countably many slices. Slices
$S_k$ are smooth manifolds diffeomorphic to $S^n$ and they converge
to the limiting slice
$S_\infty$. The set $\widetilde{X}$ is compact, but not connected.
To make it connected we define $X$ by adding to $\widetilde{X}$
a curve that connects all the sets in the family
and has the property that no part of the curve is rectifiable.

Using absolute continuity of Sobolev mappings on lines one can easily prove
the following fact.
\begin{lemma}
\label{T8}
For each $f\in W^{1,n}(S^n, X)$ we can choose a representative
(in the class of functions equal a.e.) such that
$f(S^n)\subset S_k$ for some $k=1,2,\ldots$ or $k=\infty$.
\end{lemma}
This implies the following result; for a detailed proof, see
\cite{hajlaszGAFA}.
\begin{lemma}
\label{T9}
Lipschitz mappings $\lip(S^n,X)$ are dense in $W^{1,n}(S^n,X)$.
Moreover there is $f\in W^{1,n}(S^n,X)$ and $\eps>0$ such that
if $g\in\lip(S^n,X)$, $\Vert f-g\Vert_{1,n}<\eps$, then
$f(x)\neq g(x)$ for all $x\in S^n$.
\end{lemma}
The idea of the proof is as follows. If $f(S^n)\subset S_k$ for a finite
$k$, then $f$ can be approximated by smooth mappings
$f_i\in C^\infty(S^n,S_k)$. This is known and
follows from the fact that $S_k$ is a smooth manifold.
If $f(S^n)\subset S_\infty$, then we can ``push'' the mapping a little bit,
to obtain a mapping $f_k\in W^{1,n}(S^n,S_k)$. Since the manifolds
$S_k$ converge to $S_\infty$, the mappings $f_k$ converge to $f$ in the
Sobolev norm. Now each mapping $f_k$ can be approximated by mappings in
$C^{\infty}(S^n,S_k)$ (as explained earlier)
and hence we obtain not only Lipschitz,
but even a smooth approximation of $f$.

Let $f:S^n\to S_\infty\subset X$ be a $W^{1,n}$ homeomorphism and let $\eps>0$
be as in (\ref{far}). If $g\in \lip(S^n,X)$ is such that
$\Vert f-g\Vert_{1,n}<\eps$, then $g$ cannot be a mapping into $S_\infty$,
so it must be a mapping into $S_k$ for some finite $k$ and hence
$f(x)\neq g(x)$ for all $x\in S^n$.

Now we are ready to complete the proof of Theorem~\ref{T2}.

Let $\kappa:X\to\ell^\infty$ be the Kuratowski embedding. Let
$f:S^n\to S_\infty\subset X$ be a $W^{1,n}$ homeomorphism.
We can, of course, assume that $\nabla f\neq 0$ a.e. Then
$\bar{f}=\kappa\circ f\in W^{1,n}_\kappa(S^n,X)$.
Suppose that $\bar{g}_k\in\lip(S^n,\kappa(X))$ converge to $\bar{f}$
in the Sobolev norm, $\Vert \bar{f}-\bar{g}_k\Vert_{1,n}\to 0$.
Define $g_k\in\lip(S^n,X)$ by
$g_k=\kappa^{-1}\circ \bar{g}_k$. It follows from Corollary~\ref{T7}
that
$$
|\nabla (\bar{f}-\bar{g}_k)|\geq
C\left(|\nabla f|+|\nabla g_k|\right)\chi_{\{f\neq g_k\}}\, .
$$
(Actually one needs a slight modification of the corollary, because now we consider
the Kuratowski embedding of $X$ and not of the entire ambient space $\bbbr^{n+2}$
of which $X$ is a subset, but the argument remains the same -- we modify
the proof of Lemma~\ref{T6} by choosing points $x_i$ and $x_{i'}$ from $X$;
we leave details to the reader.)
Hence
$$
0\leftarrow \Vert \bar{f}-\bar{g}_k\Vert_{1,n}^n
\geq
C\left( \int_{\{f\neq g_k\}} |\nabla f|^n +
\int_{\{f\neq g_k\}}|\nabla g_k|^n\right)\, .
$$
Thus $|\{ f\neq g_k\}|\to 0$ and
$$
\int_{\{f\neq g_k\}} |\nabla g_k|^n\to 0\, .
$$
This, in turn, implies that $g_k\to f$ in
$W^{1,n}(S^n,X)$. Therefore Lemma~\ref{T9} implies that for all
sufficiently large $k$,
$g_k\neq f$ everywhere,
which contradicts the fact that $|\{ f\neq g_k\}|\to 0$.
The proof is complete.
\hfill $\Box$

\end{document}